\documentclass[12pt]{article}%
\usepackage{amsfonts, amsmath, amscd}
\usepackage[psamsfonts]{amssymb}

\usepackage{amssymb}

\newcommand{\pr}{\par {\bf Proof.} }   %%%%%%%%%%%%%
\newcounter{df}%[section]                      %
\newcounter{pro}%[section]                      %
\newenvironment{pro}{\par%                     %
\refstepcounter{pro}%                          % оформление предложений
{\bf Proposition \arabic{pro}.} }{}% \arabic{section}.        %
\newcounter{rem}%[section]                     %
\newcounter{exa}%[section]                     %
\newcounter{teo}%[section]                      %
\newenvironment{teo}{\par%                     %
\refstepcounter{teo}%                          % оформление теорем
{\bf Theorem \arabic{teo}.} \it }{}% \arabic{section}.       %
\newcounter{cor}%[section]                     %
\newenvironment{cor}{\par%                    %
\refstepcounter{cor}%                         % оформление следствий
{\bf Corollary \arabic{cor}.} \it }{}%  \arabic{section}.   %
\newcounter{st}%[section]                      %
\newcounter{lem}%[section]                     %
\newenvironment{lem}{\par%                    %
\refstepcounter{lem}%                         % оформление лемм
{\bf Lemma \arabic{lem}.} \it }{}%  \arabic{section}.       %
\makeatletter                                                   %
\renewcommand{\section}{\@startsection{section}{1}%             % оформление
{\parindent}{3.5ex plus 1ex minus 0.2ex}{2.3ex plus 0.2ex}{\bf}}% заголовков
\makeatother                                                    %
\begin{document}%

\author{S.S. Gabriyelyan\footnote{The author was partially supported
 by Israel Ministry of Immigrant Absorption}}
\title{Minimally almost periodic group topology on countable torsion Abelian groups}
\date{}

\makeatletter
\renewcommand{\@makefnmark}{}
\renewcommand{\@makefntext}[1]{\parindent=1em #1}
\makeatother

\maketitle

\vspace{-\baselineskip}

%\subjclass[2000]{Primary 22A10, 43A40; Secondary 54H11}

\begin{abstract}
For any countable torsion subgroup $H$ of an unbounded Abelian group $G$ there is a complete  Hausdorff group topology $\tau$ such that $H$ is the von Neumann radical of $(G,\tau)$. In particular, any unbounded torsion countable Abelian group admits a complete Hausdorff minimally almost periodic (MinAP) group topology. A bounded  infinite Abelian group  admits a MinAP group topology if and only if all its leading Ulm-Kaplansky invariants are infinite. If, in addition,  $G$ is countably infinite , a MinAP group topology can be chosen to be complete.
\end{abstract}

\maketitle

\section{Introduction}

We shall write our Abelian groups additively. For a topological group $X$, $X^{\wedge}$ denotes the group of all continuous characters on $X$ endowed with the compact-open topology.  Denote by $\mathbf{n}(X) = \cap_{\chi\in X^{\wedge}} {\rm ker} \chi$ the von Neumann radical of $X$. If $\mathbf{n}(X) = X$, the group $X$ is called minimally almost periodic (MinAP). $X$ is called maximally almost periodic if $\mathbf{n}(X) = 0$. Let $H$ be a subgroup of $X$. The annihilator of $H$  we denote by $H^{\perp}$. $H$ is called dually closed in $X$ if for every $x\in X\setminus H$ there exists a character $\chi\in H^{\perp}$ such that $(\chi,x)\not= 1$. $H$ is named dually embedded in $X$ if every character of $H$ can be extended to a character of $X$.

A group $G$ with the discrete topology is denoted by $G_{\mathcal{D}}$. A group $G$ is called a group of infinite exponent or unbounded if $\exp G =\infty$ (where $\exp G$=the least common multiple of the orders of the elements of $G$), otherwise, i.e. $\exp G<\infty$, it is called a bounded group. Note that a bounded group $G$ has the form $G =\bigoplus_{p\in M} \bigoplus_{i=1}^{n_p} \mathbb{Z}(p^i)^{(k_{i,p})}$, where $M$ is a finite set of prime numbers. Leading Ulm-Kaplansky invariants of $G$ are the cardinal numbers $k_{n_p ,p}, p\in M$. The order of an element $g\in G$ we denote by $o(g)$. The subgroup of $G$ generated by an element $g$ is denoted by $\langle g\rangle$.

Let $X$ be a topological group and $\mathbf{u} =\{ u_n \}$ a sequence of elements of $X^{\wedge}$. Following D.~Dikranjan et al. \cite{DMT}, we denote by $s_{\mathbf{u}} (X)$ the set of all $x\in X$ such that $(u_n , x)\to 1$. Let $G$ be a subgroup of $X$. If $G=s_{\mathbf{u}} (X)$ we say that $\mathbf{u}$ {\it characterizes} $G$ and that $G$ is {\it characterized} (by $\mathbf{u}$) \cite{DMT}.

Following E.G.Zelenyuk and I.V.Protasov \cite{ZP1}, \cite{ZP2}, we say that a sequence $\mathbf{u} =\{ u_n \}$ in a group $G$ is a $T$-{\it sequence} if there is a Hausdorff group topology on $G$ for which $u_n $ converges to zero. The group $G$ equipped with the finest group topology with this property is denoted by $(G, \mathbf{u})$. We note also that, by Theorem 2.3.11 \cite{ZP2}, the group $(G, \mathbf{u})$ is complete.
Following \cite{BDM}, we say that a sequence $\mathbf{u} =\{ u_n\}$ is a $TB$-{\it sequence} in a group $G$ if there is a precompact Hausdorff group topology on $G$ in which  $u_n \to 0$. As it was noted in \cite{Ga2}, a sequence $\mathbf{u}$ is a $TB$-sequence if and only if it is a $T$-sequence and $(G, \mathbf{u})$ is maximally almost periodic.

Let $G$ be an infinite Abelian group. Denote  by $\mathcal{NR}(G)$ [respectively $\mathcal{NRC}(G)]$ the set of all subgroups $H$ of $G$  for which there exists a [respectively Complete] non-discrete Hausdorff group topology $\tau$ on $G$ such that von Neumann Radical of $(G,\tau)$ is $H$, i.e. $\mathbf{n} (G,\tau)=H$. It is clear that $\mathcal{NRC}(G)\subset \mathcal{NR}(G)$.

One of the most important properties of an Abelian topological group $G$ is the richness of its dual group $G^\wedge$ which is characterized by the von Neumann radical.  Thus, the general question of describing  the sets $\mathcal{NR}(G)$ and $\mathcal{NRC}(G)$, which is raised in \cite{Ga2}, is important. Since every infinite Abelian group admits a $TB$-sequence \cite{CRT}, then every infinite Abelian group admits a  complete non-trivial Hausdorff group topology with  trivial von Neumann radical. Thus, the trivial group $\{ 0\}$ belongs to $\mathcal{NRC}(G)$ for every infinite Abelian group $G$, i.e. $\mathcal{NRC}(G)$ is always not empty. A much deeper question is whether any infinite Abelian group admits a Hausdorff group topology with  {\it non-zero} von Neumann radical. The positive answer was given by M.Ajtai, I.Havas and J.Koml\'{o}s \cite{AHK}. Using the method of $T$-sequences, E.G.Zelenyuk and I.V.Protasov \cite{ZP1} proved that every infinite Abelian group admits a {\it complete}  Hausdorff group topology for which characters do not separate points, i.e. $\mathcal{NRC}(G)\not= \{ \{ 0\}\}$ for every infinite $G$.

It maybe the most interesting question is to describe all infinite Abelian groups $G$ such that $G\in \mathcal{NR}(G)$ (or $G\in \mathcal{NRC}(G)$). A simple example of a bounded group $G$ which does not admit any Hausdorff group topology $\tau$ such that $(G,\tau)$ is minimally almost periodic is given by D.~Remus \cite{Com}. This justifies the following problems:

\begin{itemize}
\item[] {\bf Question 1.} (Comfort's Problem 521 \cite{Com}) {\it Does every Abelian group which is not of bounded order admit a minimally almost periodic topological group topology? What about
    the countable case?}
\end{itemize}

Denote by $F$ the maximal torsion subgroup of an Abelian group $G$. Then the factor group $G/F$ is torsion free. Thus the following two special cases of Question 1 are natural.
\begin{itemize}
\item[] {\bf Question 2.} (Question 2.6.1 \cite{ZP2}) {\it Let $G$ be a torsion free group. Does there exist a Hausdorff group topology on $G$ with only zero character?}
\item[] {\bf Question 3.}  {\it Let $G$ be a torsion group. Under which conditions the group $G$ admits a Hausdorff minimally almost periodic group topology?}
\end{itemize}

Note that even the following  more particular question is also non-trivial: which Abelian groups $G$ admit a Hausdorff group topology $\tau$ such that $\mathbf{n}(G,\tau)$ is non-trivial and finite? This question  was raised by G.~Luk\'{a}cs \cite{Luk} and he called such groups almost maximally almost-periodic (AMAP). He proved that infinite direct sums and the Pr\"{u}fer group $\mathbb{Z} (p^{\infty}),$ for every prime $p\not= 2$, are AMAP. These results were generalized in \cite{Ngu}.
Obviously, an AMAP group contains a non-trivial finite subgroups, so it cannot be torsion free. Therefore, an Abelian group admitting an AMAP group topology is not torsion free. It turns out that the converse assertion is also true:
\begin{itemize}
\item[] {\bf Theorem A.} \cite{Ga2} {\it An Abelian group admits an AMAP group topology iff it is not torsion free.}
\end{itemize}

A complete characterization of those finitely generated subgroups of an infinite Abelian group $G$ which are the von Neumann radical for some Hausdorff group topology on $G$ is given in \cite{Ga3}.

The main goal of the paper is to describe all countable (finite or not) torsion subgroups of an infinite Abelian group $G$  which are contained in $\mathcal{NR}(G)$. In particular, we give a complete answer to Question 3 when $G$ is a countably infinite Abelian group. We prove the following theorems.

\begin{teo} \label{t1}
For any countable (finite or not) torsion subgroup $H$ of an infinite unbounded Abelian group $G$ there exists a complete Hausdorff group topology such that $H=\mathbf{n} (G)$, i.e. if $\exp (G) = \infty$, then $H\in \mathcal{NRC} (G)$.
\end{teo}

As a consequence we obtain the following.

\begin{cor} \label{c1}
Any countably infinite unbounded torsion Abelian group admits a complete Hausdorff MinAP group topology.
\end{cor}

\begin{teo} \label{t2}
Let $G$ be an infinite Abelian bounded group and $H$ its countable (finite or not) subgroup. Then the following assertions are equivalent.
\begin{enumerate}
\item Let $N=\exp (H)$. Then the group $G$ contains a subgroup of the form $\mathbb{Z} (N)^{(\omega)}$.
\item $H\in \mathcal{NRC}(G)$;
\item $H\in \mathcal{NR}(G)$.
\end{enumerate}
\end{teo}

The following corollaries are evident.

\begin{cor} \label{c2}
Let $G$ be a countably infinite bounded Abelian group $G$. Then $G$ admits a Hausdorff minimally almost periodic group topology if and only if all its leading Ulm-Kaplansky invariants are infinite. In such a case, a MinAP group topology can be chosen to be complete.
\end{cor}

\begin{cor} \label{c3}
All countable (finite or not) subgroups of an infinite bounded Abelian group $G$ belong to $\mathcal{NR} (G)$ if and only if all leading Ulm-Kaplansky invariants of $G$ are infinite. In such a case, all countable (finite or not) subgroups of $G$ belong to $\mathcal{NRC} (G)$.
\end{cor}

Note that Theorems \ref{t1} and \ref{t2} for a finite $H$ are proved in \cite{Ga3}. We will prove Theorems \ref{t1} and \ref{t2} using the following assertions.

\begin{pro} \label{p1}
{\it Let $G=\mathbb{Z}\oplus H$, where $H$ be a countable Abelian group. Then $G$ admits two complete Hausdorff group topologies $\tau_1$ and $\tau_2$ such that $\mathbf{n} (G,\tau_1) = H$ and $\mathbf{n} (G,\tau_2) = G$.}
\end{pro}

As a trivial consequence of this proposition and the Reduction Principle (see below) we obtain the next assertion.
\begin{cor} \label{c4}
Let $H$ be a countable subgroup of an infinite Abelian group $G$. If there is an element $e_0$ of infinite order in $G$ such that $\langle e_0 \rangle \cap H= \{ 0\}$, then $H\in \mathcal{NRC} (G)$.
\end{cor}

\begin{pro} \label{p2}
{\it Let $G=\mathbb{Z} (p^\infty) \oplus H$, where $H$ be a countable Abelian group. Then for any subgroup $H_1$ of $G$ containing $H$ there is a complete Hausdorff group topology $\tau$  such that $\mathbf{n} (G,\tau) = H_1$.}
\end{pro}

Since, by Baer's Theorem (Theorem 18.1 \cite{Fuc}), groups of the form $\mathbb{Z} (p^\infty)$ are direct summands of any Abelian group, we obtain the following.
\begin{cor} \label{c5}
Let $H$ be a countable subgroup of an infinite Abelian group $G$. If $G$ contains a nonzero element of infinite height, then $H\in \mathcal{NRC} (G)$.
\end{cor}

The next assertion is a generalization of Example 2.6.2 \cite{ZP2}.
\begin{pro} \label{p3}
{\it Let $G= \bigoplus_{j=0}^\infty \langle e_j \rangle$. Assume that $u_j := o(e_j) <\infty$ for every $j\geqslant 0$. Let $e'_j \in \langle e_j \rangle$ be arbitrary for every $0\leqslant j <\infty$. Set $H= \bigoplus_{j=0}^\infty \langle e'_j \rangle$. If either $a)$ $u_{j_0} =u_{j_0 +1} =\dots$, for some $j_0 \geqslant 0$ and $u_{j_0}$ is divided on every $u_0,\dots, u_{j_0 -1}$, or $b)$ $\sup\{ u_j , j\geqslant 0\} = \infty$, then $G$ admits a complete Hausdorff group topology $\tau$ such that $\mathbf{n}(G,\tau) = H$.}
\end{pro}

It turns out that the conclusion of Corollary \ref{c2} remains true also for non countable groups.
\begin{teo} \label{t3}
An infinite bounded Abelian group admits a MinAP group topology if and only if all its leading Ulm-Kaplansky invariants are infinite.
\end{teo}

We prove Theorem \ref{t3}  using the next proposition.
\begin{pro} \label{p4}
{\it Let $G= \bigoplus_{\alpha \in I} G_\alpha$, where $I$ be arbitrary set of indexes and each group
$G_\alpha$ admit a MinAP group topology. Then $G$ also admits a MinAP group topology.}
\end{pro}

\section{The proofs}

The plan of the proofs is as follows.

{\bf 1. The reduction to countable groups.}
The following lemma plays an important role.
\begin{lem} \label{l1}
{\rm \cite{Ga1}} {\it Let $K$ be a dually closed and dually embedded subgroup of a topological group $G$. Then $\mathbf{n} (K) =\mathbf{n} (G)$. }
\end{lem}

{\bf Reduction Principle}. Let $H$ be a subgroup of an infinite Abelian group $G$ and let $Y$ be arbitrary subgroup of $G$ containing $H$. If $Y$ admits a Hausdorff group topology $\tau$ such that $\mathbf{n} (Y,\tau) =H$, then we can consider the topology on $G$ in which $Y$ is open. Since any open subgroup is dually closed and dually embedded (Lemma 3.3 \cite{Nob}), by Lemma \ref{l1}, we obtain that $\mathbf{n} (G)=H$ too. Thus we can reduce our consideration to $Y$ only.

Since $H$ is countable, by this reduction principle we can split the proof of Theorem \ref{t1} in the following cases:
\begin{enumerate}
\item[(1)] {\it $H$ is unbounded and contains  nonzero elements of infinite height.} Then $H$ has the form $H= \mathbb{Z} (p^\infty) \oplus H_1$ for some prime $p$ and a countable (finite or not) subgroup $H_1$. Set $G=H$. Then Theorem \ref{t1} follows from Proposition \ref{p2} immediately.
\item[(2)] {\it $H$ is unbounded and does not contain nonzero elements of infinite height.} By the  Pr\"{u}fer Theorem (Theorem 11.3 \cite{Fuc}), $H=\bigoplus_{j=0}^\infty \langle e_j \rangle,$
     where $u_j := o(e_j) <\infty$ for every $j\geqslant 0$ and $\sup\{ u_j , j\geqslant 0\} = \infty$. Then Theorem \ref{t1} immediately follows from Proposition \ref{p3}.
\item[(3)] {\it $H$ is bounded and $G$ contains an element of infinite order $e_0$.} Since $\langle e_0\rangle \cap H=\{ 0\}$, we can put $G= \mathbb{Z} \oplus H$. So Theorem \ref{t1} follows from Proposition \ref{p1}.
\item[(4)] {\it $H$ is bounded and $G$ contains a nonzero element of infinite height.} Thus $G$ contains a subgroup of the form $\mathbb{Z} (p^\infty)$ for some prime $p$. Hence $\mathbb{Z} (p^\infty) \cap H$ is either zero or $H=\mathbb{Z} (p^k) \oplus H_1$, where $\mathbb{Z} (p^k)\subset \mathbb{Z} (p^\infty)$ and $\mathbb{Z} (p^\infty) \cap H_1 = \{ 0\}$. In the first case  we put $G= \mathbb{Z} (p^\infty) \oplus H$ and in the second one we set $G= \mathbb{Z} (p^\infty) \oplus H_1$. Then Theorem \ref{t2} follows directly from Proposition \ref{p2}.
\item[(5)] {\it $H$ is bounded, $G$ is unbounded and contains no nonzero elements of infinite height and infinite order.}  By the  Pr\"{u}fer Theorem (Theorem 11.2 \cite{Fuc}), there is a finite set $M$ of prime numbers such that
    \[
    H =\bigoplus_{p\in M} \bigoplus_{i=1}^{n_p} \mathbb{Z}(p^{a_{i,p}})^{(k_{i,p})}, \mbox{ where } 1\leqslant a_{1,p} <\dots < a_{n_p,p}, \; 1\leqslant k_{i,p} \leqslant\infty, \; p\in M.
    \]
For every $0\leqslant j < k_{i,p}$ let $e'_{j,i,p}$ be a generator of $\mathbb{Z}(p^{a_{i,p}})$. Then
\[
H =\bigoplus_{p\in M} \bigoplus_{i=1}^{n_p} \bigoplus_{0\leqslant j< k_{i,p}} \langle e'_{j,i,p} \rangle .
\]
Let us denote by $b_{j,i,p}$ the height of $e'_{j,i,p}$ in the $p$-component $G_p$ of $G$. Then  we have the following two cases.

 a) $\sup\{ b_{j,i,p}\} =\infty$. Denote by $e_{j,i,p}$ an element of $G_p$ such that $p^{b_{j,i,p}} e_{j,i,p} = e'_{j,i,p}$. Then all $e_{j,i,p}$ are independent. Set
\[
G_1 = \bigoplus_{p\in M} \bigoplus_{i=1}^{n_p} \bigoplus_{0\leqslant j< k_{i,p}} \langle e_{j,i,p} \rangle .
\]
Then Theorem \ref{t1} follows from Proposition \ref{p3}(b).

b) $\sup\{ b_{j,i,p}\} <\infty$. Since $G$ is unbounded, there is a sequence $\{ e_j\}$ of elements in $G$ such that $o(e_j)\to \infty$ and all $e'_{j,i,p}$ and $e_j$ are independent. Then  $G$ contains the following  direct sum of the groups
\[
G_1 = H\oplus \bigoplus_{j=0}^\infty \langle e_j \rangle = \bigoplus_{p\in M} \bigoplus_{i=1}^{n_p} \bigoplus_{0\leqslant j< k_{i,p}} \langle e'_{j,i,p} \rangle \oplus \bigoplus_{j=0}^\infty \langle e_j \rangle .
\]
So also in this case Theorem \ref{t1} follows from Proposition \ref{p3}(b).
\end{enumerate}

Thus to prove Theorem \ref{t1} we need only to prove Propositions \ref{p1}-\ref{p3}.

{\bf 2. Construction of a topology.} We construct a special $T$-sequence $\mathbf{d}$ on $G$ and equip $G$ under the Protasov-Zelenyuk group topology generated by the choosing $T$-sequence $\mathbf{d}$.

{\bf 3. Computation of the von Neumann kernel.} Using the following theorem
\begin{itemize}
\item[] {\bf Theorem B.} \cite{Ga1} {\it If $\mathbf{d}=\{ d_n\}$ is a $T$-sequence of an Abelian group $G$, then $\mathbf{n} (G,\mathbf{d})= s_{\mathbf{d}} \left( (G_{\mathcal{D}})^\wedge \right)^{\perp}$ (algebraically).}
\end{itemize}
we prove that $s_{\mathbf{d}} \left( (G_{\mathcal{D}})^\wedge \right)^\perp =H$. $\Box$

In the sequel we need some notations. For a sequence $\{ d_n \}$ and $l,m\in \mathbb{N}$, one puts \cite{ZP1}:
\[
\begin{split}
A(k,m) = & \left\{ n_1 d_{r_1} +\dots +n_s d_{r_s} | \; \mbox{where } \right. \\
 &  \left. m\leqslant r_1 <\dots < r_s , n_i \in \mathbb{Z} \setminus \{ 0\} , \sum_{i=1}^s | n_i | \leqslant k+1 \right\} \cup \{ 0\}.
\end{split}
\]

Denote by $S_n =1+\dots +n=\frac{n(n+1)}{2}$ and for every $n\geqslant 1$ we put
$t(n) := \max\{ t:\; n\geqslant S_t\}$.

For a prime $p$ and $n\in \mathbb{N}$ we set
\[
f_n =p^{n^3 - n^2} +\dots + p^{n^3 - 2n} +  p^{n^3 - n} + p^{n^3} \in \mathbb{Z},
\]
then $f_n < 2 p^{n^3} \leqslant p^{n^3 +1}$. For any $0<r_1 <r_2 <\dots <r_v$ and integers $l_1, l_2,\dots , l_v$ such that $\sum_{i=1}^v |l_i | \leqslant k+1$, we have
\begin{equation} \label{1}
| l_1 f_{r_1 } + l_2 f_{r_2 } +\dots + l_v f_{r_v } | < (k+1) f_{r_v} \leqslant (k+1) p^{r_v^3 +1} .
\end{equation}

The following lemma is a  generalization of Lemma 2.2 in \cite{Ga3}.

\begin{lem} \label{l2}
Let $G=\mathbb{Z} \oplus H$, where $H$ be a countable Abelian group, and $e_0$ be the generator of $\mathbb{Z}$. Let $H\setminus \{ 0\} = \{ e_n \}_{n=1}^{|H|-1}$ be any enumeration of the non-zero elements of $H$. Choose  a prime $p$ and arbitrary sequence $\varepsilon_n$ such that $\varepsilon_n \in \{ -1,0,1\}$ for all $n\geqslant 3$. Define the following sequence ($n\geqslant 3$)
\[
\begin{split}
\mbox{ if } |H|=\infty : & \quad d_{2n} =  p^n e_0 , \;
d_{2n-1} =  f_n e_0   + \varepsilon_n e_{n\left({\rm mod}\, S_{t(n)} \right) } ,  \\
\mbox{ if } |H|<\infty : & \quad d_{2n} =  p^n e_0 , \;
d_{2n-1} =  f_n e_0   + \varepsilon_n  e_{n\left({\rm mod} |H| \right) } .
\end{split}
\]
Then the sequence $\mathbf{d}=\{ d_n \}$ is a $T$-sequence.
\end{lem}

\pr Let $0\not= g= b e_0+ \epsilon e_{q} \in G$,  where $b\in \mathbb{Z}$ and $\epsilon$ is either $0$ or $1$. Let $k\geqslant 0$. Set $m=20(|b| +1) (k+1)$. By Theorem 2.1.4 \cite{ZP2}, it is enough to prove that $g\not\in A(k,m)$.

Denote by $\eta_n = S_{t(n)}$ if $ |H|=\infty$ and $\eta_n = |H|$ if $ |H|<\infty$ for every $n\geqslant 3$.

a) Let $\sigma \in A(k,m)$ have the following form
\[
\sigma =l_1 d_{2r_1 } + l_2 d_{2r_2 } +\dots + l_s d_{2r_s} = (l_1 p^{r_1}  +\dots + l_s p^{r_s}) e_0 = p^{r_1} \cdot\sigma' \cdot e_0,
\]
where  $m\leqslant 2r_1 <2r_2 <\dots <2r_s$  and  $\sigma' \in \mathbb{Z}$.
If $\sigma' =0$, then $\sigma \not= g$. If $\sigma' \not= 0$, then $p^{r_1} > p^{10|b|} >|b|$, and $\sigma \not= g.$

b) Let $\sigma \in A(k,m)$ have the following form
\[
\sigma =l_1 d_{2r_1 -1} + l_2 d_{2r_2 -1} +\dots + l_s d_{2r_s -1},
\]
where $m< 2r_1 -1 < 2r_2 -1 < \dots < 2r_s -1 $ and nonzero integers $l_1, l_2,\dots , l_s$ be such that   $\sum_{i=1}^s |l_i | \leqslant k+1$. Then
\[
\sigma = (l_1 f_{r_1 } + \dots + l_{s-1} f_{r_{s-1} } + l_s f_{r_s } )e_0 +
    l_1 \varepsilon_{r_1} e_{r_1\left({\rm mod}\, \eta_{r_1} \right) } +\dots + l_s \varepsilon_{r_s} e_{r_s\left({\rm mod}\, \eta_{r_s} \right) }.
\]
Since $n^3 < (n+1)^3 -(n+1)^2$ and $r_s >|b|+(k+1)$, by (\ref{1}), we can estimate the coefficient $\phi_0$ of $e_0$ in $\sigma$ as follows
\[
\begin{split}
|\phi_0|\geqslant & |l_1 f_{r_1 } + \dots + l_{s-1} f_{r_{s-1} } + l_s f_{r_s } | - (k+1)
>  f_{r_s } - k\cdot p^{r_{s-1}^3 +1} -(k+1) \\
= &
p^{r_{s}^3 } + \left( p^{r_{s}^3 -r_s }+\dots + p^{r_{s}^3 -r_s^2} - k\cdot p^{r_{s-1}^3 +1} -k-1\right) >  p^{r_{s}^3 } >|b|.
\end{split}
\]
Hence $\phi_0 \not= b$ and $\sigma \not= g$.

c) Let $\sigma \in A(k,m)$ have the following form
\[
\sigma = l_1 d_{2r_1 -1} + l_2 d_{2r_2 -1} +\dots + l_s d_{2r_s -1} +l_{s+1} d_{2r_{s+1}}  +\dots + l_{h} d_{2r_{h}},
\]
where $0<s<h$ and
\[
\begin{split}
  m< & 2r_1 -1 < 2r_2 -1 < \dots < 2r_s -1 ,   \\
  m\leqslant & 2r_{s+1} < 2r_{s+2} < \dots < 2r_{h}, \quad l_i \in \mathbb{Z} \setminus \{ 0\} , \sum_{i=1}^h |l_i | \leqslant k+1.
\end{split}
\]
Since the number of summands in $d_{2r_s -1}$ is at least $r_s +1> 10(k+1)$ and $h-s <k+1$, there exists $r_s - 2> i_0 >2$ such that for every $1\leqslant w \leqslant h-s$ we have
\[
\mbox{either } r_{s+w} < r_s^3 - (i_0 +2)r_s  \mbox{ or }  r_{s+w} > r_s^3 - (i_0 -1)r_s .
\]
The set of all $w$ such that $r_{s+w} < r_s^3 - (i_0 +2)r_s$ we denote by $B$ (it may be empty or has the form $\{ 1,\dots ,\delta \}$ for some $1\leqslant \delta \leqslant h-s$). Set $D =\{ 1,\dots ,h-s \} \setminus B$. Thus
\[
\begin{split}
  \sigma =  &  l_1 \varepsilon_{r_1} e_{r_1\left({\rm mod}\, \eta_{r_1} \right) } +\dots + l_s \varepsilon_{r_s} e_{r_s\left({\rm mod}\, \eta_{r_s} \right) } +  \left( l_1 f_{r_1 } + \dots + l_{s-1} f_{r_{s-1} }\right) e_0  \\
   & + \sum_{w\in B} l_{s+w} d_{2r_{s+w}} + \left( l_s p^{r_s^3 - r_s^2} +\dots + l_s p^{r_s^3 - (i_0 +2)r_s}\right) e_0  \\
   & + \left( l_s p^{r_s^3 - (i_0 +1)r_s} + l_s  p^{r_s^3 - i_0 r_s}\right) e_0 +
     \left[ \left( l_s p^{r_s^3 - (i_0 -1)r_s} +\dots + l_s  p^{r_s^3} \right) e_0 + \sum_{w\in D} l_{s+w} d_{2r_{s+w}}\right] .
\end{split}
\]
We can estimate the coefficient of $e_0$  in row 2, which we denote by $A_2$, as follows
\begin{equation} \label{2}
|A_2 |   < \sum_{w\in B} |l_{s+w} | p^{r_s^3 - (i_0 +2)r_s}+ |l_s | 2 p^{r_s^3 - (i_0 +2)r_s} <3(k+1) p^{r_s^3 - (i_0 +2)r_s} <p^{r_s^3 - (i_0 +1)r_s},
\end{equation}
since $3(k+1) <r_s < p^{r_s} -1$. For the coefficient of $e_0$ of the second summand in row 3, which we denote by $A_3$, we have
\begin{equation} \label{3}
A_3 = l_s p^{r_s^3 - (i_0 -1)r_s} +\dots + l_s  p^{r_s^3} + \sum_{w\in D} l_{s+w} p^{r_{s+w}} =  p^{r_s^3 - (i_0 -1)r_s} \cdot \sigma'' ,
\end{equation}
where  $\sigma'' \in \mathbb{Z}$. Set $\phi_0$ is the coefficient of $e_0$ in $\sigma$. Let $\sigma'' \not= 0$. By (\ref{1})-(\ref{3}), we can estimate  $\phi_0$ from below as follows:
\[
\begin{split}
 |\phi_0| \geqslant & -(k+1) + |l_1 f_{r_1 } + \dots + l_{s-1} f_{r_{s-1} } + A_2 + l_s p^{r_s^3 - (i_0 +1)r_s} + l_s  p^{r_s^3 - i_0 r_s} +  A_3 |  \\
 >  & p^{r_s^3 - (i_0 -1)r_s } -k\cdot p^{r_{s-1}^3 +1} -2(k+1) p^{r_s^3 - i_0 r_s} -(k+1) \\
 >  & p^{r_s^3 - (i_0 -1)r_s } -3(k+1) p^{r_s^3 - i_0 r_s} > p^{r_s^3 - i_0 r_s }>  p^{r_s^2} > p^{|b|} >|b|.
\end{split}
\]
Hence $\phi_0 \not= b$ and $\sigma \not= g$. Let $\sigma'' = 0$. By (\ref{1}) and (\ref{2}), we have
\[
\begin{split}
  |\phi_0| \geqslant & -(k+1) + |l_1 f_{r_1 } + \dots + l_{s-1} f_{r_{s-1} } +A_2  + l_s p^{r_s^3 - (i_0 +1)r_s} + l_s  p^{r_s^3 - i_0 r_s} |  \\
 >  & p^{r_s^3 - i_0 r_s} - k\cdot p^{r_{s-1}^3 +1} -(k+2) p^{r_s^3 - (i_0 +1)r_s} -(k+1)  \\
 >  & p^{r_s^3 - i_0 r_s} -3(k+1) p^{r_s^3 - (i_0 +1)r_s } > p^{r_s^3 - (i_0 +1)r_s } >p^{r_s^2} > p^{|b|} >|b|.
\end{split}
\]
Hence $\phi_0 \not= b$ and $\sigma \not= g$ too. Thus $\mathbf{d}$ is a $T$-sequence. $\Box$

{\bf Proof of Proposition \ref{p1}.}${}^1$ \footnote{${}^1$If $H$ is finite, the proposition is proved in \cite{Ga3}.}  Let $\mathbf{d} =\{ d_n\}$ be a $T$-sequence defined in Lemma \ref{l2}. By Theorem B, it is enough to prove that $s_\mathbf{d} \left( (G_{\mathcal{D}})^\wedge \right) = \{ 0\}$.

Let
\[
\omega = x+y \in (G_{\mathcal{D}})^\wedge , \mbox{ where } x \in \mathbb{T},  y\in (H_{\mathcal{D}})^\wedge  \mbox{ and } \; (d_n , \omega) \to 1.
\]
Then $(d_{2n} , \omega)= (p^n e_0 ,x)\to 1$. Hence $x \in \mathbb{Z} (p^{\infty})$ (see \cite{Arm} or Remark 3.8 \cite{BD2}). Let $x = \frac{\rho}{p^{\tau}} , \rho\in \mathbb{Z}, \tau >0$.

Denote by $\nu_n = S_{t(n)}$ if $ |H|=\infty$ and $\nu_n = n|H|$ if $ |H|<\infty$ for every $n\geqslant 3$. The sequence $\eta_n$ is defined in Lemma \ref{l2}.

1) {\it Assume that $\varepsilon_n =1$ for every $n\geqslant 3$.} Then
\[
d_{2n-1} =  f_n e_0 + e_{n\left({\rm mod}\, \eta_n \right) }
\]
and for any $n > S_\tau$ we have $(d_{2\nu_n -1} , \omega)= (e_0,x) \to 1$ and hence $x=0$. Now let $y\not= 0$. Then  we can choose $j>0$ such that $(e_j ,y)\not= 1$. Hence for any $n > S_{\max\{\tau, j\}}$ we obtain
\[
(d_{2(\nu_n +j)-1} , \omega)= (e_j ,y) \not\to 1.
\]
It is a contradiction. So also $y=0$. Thus  $s_\mathbf{d} \left( (G_{\mathcal{D}})^\wedge \right) = \{ 0\}$. By Theorem B, $\mathbf{n} (G) =G$.

2) {\it Assume that $\varepsilon_n =1$ for every $n\not= \nu_k, k\in \mathbb{N},$ and $\varepsilon_n =0$ if $n= \nu_k$ for some $k$.} Now if $n= \nu_k$ for some $k>S_\tau$, we obtain $(d_{2n-1}, \omega)= (d_{2\nu_k -1} , \omega)= (f_{\nu_k} e_0 , x)=1$. Let $n\not= \nu_k, k\in \mathbb{N}$. Then for $n>S_\tau$ we have
\[
(d_{2n-1} ,\omega) = ( f_n e_0 + e_{n\left({\rm mod}\, \eta_n \right) }, \omega) =(e_{n\left({\rm mod}\, \eta_n \right) }, y).
\]
Now let $y\not= 0$. Then  we can choose $j>0$ such that $(e_j ,y)\not= 1$. Hence for any $n > S_{\max\{\tau, j\}}$ we obtain $(d_{2(\nu_n +j)-1} , \omega)= (e_j ,y) \not\to 1$.
It is a contradiction. So  $y=0$. Thus  $s_\mathbf{d} \left( (G_{\mathcal{D}})^\wedge \right) = \mathbb{Z} (p^{\infty})\subset \mathbb{T}$. By Theorem B, $\mathbf{n} (G) =H$.
$\Box$

Putting
\[
\widetilde{f}_n =\frac{1}{p^{n^3 - n^2}} +\dots + \frac{1}{p^{n^3 - 2n}} +  \frac{1}{p^{n^3 - n}} + \frac{1}{p^{n^3}} \in \mathbb{Z} (p^{\infty}),
\]
we obtain \cite{Luk}
\begin{equation} \label{4}
0 < \widetilde{f}_n =\frac{1}{p^{n^3 - n^2}} +\dots + \frac{1}{p^{n^3 - 2n}} +  \frac{1}{p^{n^3 - n}} +  \frac{1}{p^{n^3}} < \frac{n+1}{p^{n^3 - n^2}} \to 0,
\end{equation}
By (\ref{4}), we have
\[
\begin{split}
{} & 0< \widetilde{f}_2 + \widetilde{f}_4 +\dots <\frac{3}{p^4} + \frac{5}{p^{48}} +\dots <1, \\
{} & \widetilde{f}_{2n+2} + \widetilde{f}_{2n+4} +\dots < \frac{2n+3}{p^{(2n+2)^3 - (2n+2)^2}} + \frac{2n+5}{p^{(2n+4)^3 - (2n+4)^2}} +\dots < \frac{1}{p^{(2n)^3 - (2n)^2}} < \widetilde{f}_{2n} .
\end{split}
\]
So for any distinct sequences $0< j'_1 < j'_2 <\dots$ and $0< j''_1 < j''_2 <\dots$ the numbers $\widetilde{f}_{2j'_1} + \widetilde{f}_{2j'_2} +\dots$ and $\widetilde{f}_{2j''_1} + \widetilde{f}_{2j''_2}  +\dots$ also are distinct. Thus we can choose a sequence $0< j_1 < j_2 <\dots$ such that the number
\[
\beta:= \widetilde{f}_{2j_1} + \widetilde{f}_{2j_2} + \widetilde{f}_{2j_3} +\dots
\]
is irrational. Now we set
\begin{equation} \label{5}
\beta_n := \widetilde{f}_{2j_1} + \widetilde{f}_{2j_2} + \dots + \widetilde{f}_{2j_n}.
\end{equation}
Thus, if we consider the group $\mathbb{Z} (p^{\infty})$ as a subgroup of $\mathbb{T}$, then the sequence $\{ \beta_n\} \subset \mathbb{Z} (p^{\infty})$ converges to $\beta$ in $\mathbb{T}$, $\langle \beta\rangle$ is dense in $\mathbb{T}$ and $\langle \beta\rangle \cap \mathbb{Z} (p^{\infty}) =\{ 0\}$.

\begin{lem} \label{l3}
Let $G=\mathbb{Z} (p^{\infty})\oplus H$, where $H$ be a countable infinite Abelian group. Let  $\{ e_i\}_{i=0}^{|H|-2}$ be any enumeration of the set $H\setminus \{ 0\}$. Define the following sequence ($n\geqslant 3$)
\[
\begin{split}
\mbox{ if } |H|=\infty : & \quad d_{3n} =  \frac{1}{p^{n}} \in \mathbb{Z} (p^{\infty}), \;
d_{3n+1} =  \widetilde{f}_{2j_n +1}  + e_{n\left({\rm mod}\, S_{t(n)} \right) } , \;
d_{3n+2} =  \beta_n , \\
\mbox{ if } |H|<\infty : & \quad d_{3n} =  \frac{1}{p^{n}} \in \mathbb{Z} (p^{\infty}), \;
d_{3n+1} =  \widetilde{f}_{2j_n +1}  + e_{n\left({\rm mod} |H|-1 \right) } , \;
d_{3n+2} =  \beta_n .
\end{split}
\]
where $\widetilde{f}_n \in \mathbb{Z} (p^{\infty})$ and $\beta_n \in \mathbb{Z} (p^{\infty})$ is defined in (\ref{5}). Then $\mathbf{d} =\{ d_n\}$ is a $T$-sequence.
\end{lem}

\pr Let $g= \frac{b}{p^z} + \varepsilon e_{j} \not= 0$, where  $\varepsilon\in \{ 0,1\}$ and either $0\not= b\in \mathbb{Z}$ and $\frac{b}{p^z} \in \mathbb{Z} (p^{\infty})$ be an irreducible fraction or $b=0$. Let $k\geqslant 0$.  By Theorem 2.1.4 \cite{ZP2}, it is enough to prove that $g\not\in A(k,m)$ for some $m$. For $\sigma \in G$ we denote by $\phi_0$ the element of $\mathbb{Z} (p^{\infty})$ such that $\sigma - \phi_0 \in H$.

Set $m= 30p(k+1)(z+1)$. We will prove that $m$ is desired.

a) Let $0\not= \sigma \in A(k,m)$ have the following form
\[
\sigma =l_1 d_{3r_1 } + l_2 d_{3r_2 } +\dots + l_s d_{3r_s} ,
\]
where  $m\leqslant 3r_1 <3r_2 <\dots <3r_s$. If $b=0$, then $\sigma \not= g$. If $b\not= 0$, then
\[
0< |\phi_0 | =  |l_{1} d_{3r_{1}} +l_{2} d_{3r_{2}} +\dots + l_{s} d_{3r_{s}} |
\leqslant  \sum_{i=1}^s \frac{|l_i|}{p^{r_i}} \leqslant  \frac{k+1}{p^{r_1}} < \frac{k+1}{p^{k+1+z}} < \frac{1}{p^{z}} \leqslant \frac{|b|}{p^z}.
\]
So $\phi_0\not= \frac{b}{p^z}$ and $\sigma \not= g$.

b) Let $\sigma \in A(k,m)$ have the following form
\[
\sigma = l_1 d_{3r_1 +1} + \dots + l_\mu d_{3r_\mu +1} +l_{\mu +1} d_{3r_{\mu+1} +2}  +\dots + l_{s} d_{3r_{s} +2},
\]
where $0\leqslant \mu \leqslant s$ and
\[
\begin{split}
 {} & m< 3r_1 +1 < 3r_2 +1 < \dots < 3r_\mu +1 ,  \mbox{ if } \mu\not=0, \\
{} &  m< 3r_{\mu +1}+2 < 3r_{\mu +2} +2 < \dots < 3r_{s}+2 ,  \mbox{ if } \mu < s, \; l_i \in \mathbb{Z} \setminus \{ 0\} ,  \sum_{i=1}^h |l_i | \leqslant k+1 .
\end{split}
\]

Set $\nu = 2j_{r_\mu} +1$ and $l_\nu = l_\mu$ if $r_\mu \geqslant r_s$ and $\nu = 2j_{r_s}$ and $l_\nu = l_s$ otherwise.
Since $n^3 < (n+1)^3 -(n+1)^2$ and $\min \{ r_\mu , r_s\} >5p(k+1)$, we have
\[
\phi_0 = \frac{z'}{p^{\nu^3 - \nu}} + \frac{l_\nu}{p^{\nu^3 }}, \mbox{ where } z' \in \mathbb{Z}.
\]
Since $|l_\nu | \leqslant k+1< \frac{\nu}{p} < p^{\nu -1}$, then $\phi_0 \not= 0$ and we have the following: if $\phi_0 =\frac{z''}{p^{\alpha}}, z'' \in \mathbb{Z}$, is an irreducible fraction, then $\alpha > \nu^3 - \nu +1 > 5z$. Hence  $\phi_0\not= \frac{b}{p^z}$ and $\sigma \not= g$.

c) Let $\sigma \in A(k,m)$ have the following form
\[
\sigma = l_1 d_{3r_1 +1} + \dots + l_\mu d_{3r_\mu +1} +l_{\mu +1} d_{3r_{\mu+1} +2}  +\dots + l_{s} d_{3r_{s} +2} +l_{s+1} d_{3r_{s+1}}  +\dots + l_{h} d_{3r_{h}},
\]
where $0< s<h$ and
\[
\begin{split}
{} & m< 3r_1 +1 < 3r_2 +1 < \dots < 3r_\mu +1 ,  \mbox{ if } \mu\not=0, \\
{} &  m< 3r_{\mu +1}+2 < 3r_{\mu +2} +2 < \dots < 3r_{s}+2 ,  \mbox{ if } \mu < s, \\
{} &  m\leqslant 3r_{s+1} < 3r_{s+2} < \dots < 3r_{h}, \quad l_i \in \mathbb{Z} \setminus \{ 0\} ,  \sum_{i=1}^h |l_i | \leqslant k+1 .
\end{split}
\]
Also we set $\nu = 2j_{r_\mu} +1$ and $l_\nu = l_\mu$ if $r_\mu \geqslant r_s$ and $\nu = 2j_{r_s}$ and $l_\nu = l_s$ otherwise.
Since $n^3 < (n+1)^3 -(n+1)^2$, the element $\phi_1$  of $l_1 d_{3r_1 +1} + \dots + l_\mu d_{3r_\mu +1} +l_{\mu +1} d_{3r_{\mu+1} +2}  +\dots + l_{s} d_{3r_{s} +2}$, which is contained in $\mathbb{Z} (p^\infty)$, can be represented in the form
\[
\phi_1 = \frac{z''}{p^{\nu^3 - \nu^2 -1}} + l_\nu \widetilde{f}_\nu, \mbox{ where } z'' \in \mathbb{Z}.
\]
Since the number of summands in $\widetilde{f}_{\nu}$ is  $\nu +1 > \frac{2}{3} m >10p(k+1)$ and $h-s <k+1$, there exists $\nu - 2> i_0 >2$ such that for every $1\leqslant w \leqslant h-s$ we have
\[
\mbox{either } r_{s+w} < \nu^3 - (i_0 +2)\nu  \mbox{ or }  r_{s+w} > \nu^3 - (i_0 -1)\nu .
\]
The set of all $w$ such that $r_{s+w} < \nu^3 - (i_0 +2)\nu$ we denote by $K$ (it may be empty or has the form $\{ 1,\dots ,a\}$ for some $1\leqslant a \leqslant h-s$). Set $L =\{ 1,\dots ,h-s \} \setminus K$. Thus
\[
\begin{split}
  \phi_0 = & \frac{z''}{p^{\nu^3 - \nu^2 -1}} + \sum_{w\in K} l_{s+w} d_{3r_{s+w}} + \frac{l_\nu}{p^{\nu^3 - \nu^2 }} +\dots +\frac{l_\nu}{p^{\nu^3 - (i_0 +2)\nu}} \\
   & +\left( \frac{l_\nu}{p^{\nu^3 - (i_0 +1)\nu}} +\frac{l_\nu}{p^{\nu^3 - i_0 \nu}}\right)+
    \left(\frac{l_\nu}{p^{\nu^3 - (i_0 -1)\nu}} + \dots +\frac{l_\nu}{p^{\nu^3 }} + \sum_{w\in L} l_{s+w} d_{3r_{s+w}} \right).
\end{split}
\]
So the element $\phi_2$  in row 1, can be represented in the form $\frac{c}{p^{\nu^3 - (i_0 +2)\nu}}$, for some $c \in \mathbb{Z}$.
Since $\nu > 10p(k+1)$, then $\frac{1}{1- 1/p^{\nu}} < \frac{32}{31}$ and $2k < p^{2k} <p^{\nu}$. Thus
we can estimate the second summand in row 2, which is denoted by $\phi_3$, as follows
\[
\begin{split}
|\phi_3 | = & \left| \left( \frac{l_\nu}{p^{\nu^3 - (i_0 -1)\nu}} + \dots +\frac{l_\nu}{p^{\nu^3 }} \right) + \sum_{w\in L} l_{s+w} \frac{1}{p^{r_{s+w}}} \right|  <
  \frac{|l_\nu |}{p^{\nu^3 - (i_0 -1)\nu}} \cdot \frac{1}{1-\frac{1}{p^{\nu}}} + \frac{k}{p^{\nu^3 - (i_0 -1)\nu +1} }\\
  < &  \frac{1}{p^{\nu^3 - (i_0 -1)\nu}} \left( k\frac{32}{31} + k\frac{1}{p} \right) <
\frac{2k}{p^{\nu^3 - (i_0 -1)\nu}} < \frac{1}{p^{\nu^3 - i_0 \nu}}.
\end{split}
\]
So
\[
\phi_0 = \frac{c}{p^{\nu^3 - (i_0 +2)\nu}}  +\frac{l_\nu}{p^{\nu^3 - (i_0 +1)\nu}} +\frac{l_\nu}{p^{\nu^3 - i_0 \nu}} + \phi_3 \not= 0
\]
and  if $\phi_0 =\frac{c''}{p^{\alpha}}, c'' \in \mathbb{Z}$, is an irreducible fraction, then $\alpha \geqslant \nu^3 - (i_0 +1)\nu > 5z$. Hence $\phi_0 \not= \frac{b}{p^z}$ and $\sigma \not= g$. Thus $\{ d_n\}$ is a $T$-sequence. $\Box$

The following lemma is a  generalization of Lemma 2.3 in \cite{Ga3} with a similar proof.

\begin{lem} \label{l4}
Let $G=\mathbb{Z} (p^{\infty})\oplus H$, where $H$ be a countable  Abelian group. Let $e_0$ be an (zero or nonzero) element of $\mathbb{Z} (p^{\infty})$ and $\{ e_i\}_{i=1}^{|H|-1}$ be any enumeration of the set $H\setminus \{ 0\}$. Define the following sequence ($n\geqslant 3$)
\[
\begin{split}
\mbox{ if } |H|=\infty : & \quad d_{2n} =  \frac{1}{p^{n}} \in \mathbb{Z} (p^{\infty}),  \;
d_{2n-1} =  \widetilde{f}_n  + e_{n\left({\rm mod}\, S_{t(n)} \right) } , \mbox{ where } \widetilde{f}_n \in \mathbb{Z} (p^{\infty}),  \\
\mbox{ if } |H|<\infty : & \quad d_{2n} =  \frac{1}{p^{n}} \in \mathbb{Z} (p^{\infty}), \;
d_{2n-1} =  \widetilde{f}_n  + e_{n\left({\rm mod} |H| \right) } , \mbox{ where } \widetilde{f}_n \in \mathbb{Z} (p^{\infty}) .
\end{split}
\]
Then $\mathbf{d} =\{ d_n\}$ is a $T$-sequence.
\end{lem}

\pr Let $g= \frac{b}{p^z} + \varepsilon e_{j} \not= 0$, where  $\varepsilon\in \{ 0,1\}$ and either $0\not= b\in \mathbb{Z}$ and $\frac{b}{p^z} \in \mathbb{Z} (p^{\infty})$ be an irreducible fraction or $b=0$. Let $k\geqslant 0$.  By Theorem 2.1.4 \cite{ZP2}, it is enough to prove that $g\not\in A(k,m)$ for some $m$. For $\sigma \in G$ we denote by $\phi_0$ the element of $\mathbb{Z} (p^{\infty})$ such that $\sigma - \phi_0 \in H$.

If $e_0 \not= 0$, we assume that it has order $p^\delta$. If $e_0 = 0$, we set $\delta =0$.

Set $m= 30p(k+1)(z+1) + \delta$. We will prove that $m$ is desired.

a) Let $0\not= \sigma \in A(k,m)$ have the following form
\[
\sigma =l_1 d_{2r_1 } + l_2 d_{2r_2 } +\dots + l_s d_{2r_s} ,
\]
where  $m\leqslant 2r_1 <2r_2 <\dots <2r_s$. If $b=0$, then $\sigma \not= g$. If $b\not= 0$, then
\[
0< |\phi_0 | =  |l_{1} d_{2r_{1}} +l_{2} d_{2r_{2}} +\dots + l_{s} d_{2r_{s}} |
\leqslant  \sum_{i=1}^s \frac{|l_i|}{p^{r_i}} \leqslant  \frac{k+1}{p^{r_1}} < \frac{k+1}{p^{k+1+z}} < \frac{1}{p^{z}} \leqslant \frac{|b|}{p^z}.
\]
So $\phi_0\not= \frac{b}{p^z}$ and $\sigma \not= g$.

b) Let $\sigma \in A(k,m)$ have the following form
\[
\sigma =l_1 d_{2r_1 -1} + l_2 d_{2r_2 -1} +\dots + l_s d_{2r_s -1},
\]
where $m< 2r_1 -1 < 2r_2 -1 < \dots < 2r_s -1 $ and integers $l_1, l_2,\dots , l_s$ be such that $l_s \not= 0$ and  $\sum_{i=1}^s |l_i | \leqslant k+1$. Since $n^3 < (n+1)^3 -(n+1)^2$ and $r_s >5p(k+1)$, we have
\[
\phi_0 = \frac{z'}{p^{r_s^3 - r_s}} + \frac{l_s}{p^{r_s^3 }}, \mbox{ where } z' \in \mathbb{Z}.
\]
Since $|l_s | \leqslant k+1< \frac{r_s}{p} < p^{r_s -1}$, then $\phi_0 \not= 0$ and we have the following: if $\phi_0 =\frac{z''}{p^{\alpha}}, z'' \in \mathbb{Z}$, is an irreducible fraction, then $\alpha > r^3_{s} - r_s +1 > 5z$. Hence  $\phi_0\not= \frac{b}{p^z}$ and $\sigma \not= g$.

c) Let $\sigma \in A(k,m)$ have the following form
\[
\sigma = l_1 d_{2r_1 -1} + l_2 d_{2r_2 -1} +\dots + l_s d_{2r_s -1} +l_{s+1} d_{2r_{s+1}}  +\dots + l_{h} d_{2r_{h}},
\]
where $0< s<h$ and
\[
\begin{split}
 {} & m< 2r_1 -1 < 2r_2 -1 < \dots < 2r_s -1 ,   \\
{} &  m\leqslant 2r_{s+1} < 2r_{s+2} < \dots < 2r_{h}, \quad l_i \in \mathbb{Z} \setminus \{ 0\} ,  \sum_{i=1}^h |l_i | \leqslant k+1 .
\end{split}
\]
Since $n^3 < (n+1)^3 -(n+1)^2$, the element $\phi_1$  of $l_1 d_{2r_1 -1} + l_2 d_{2r_2 -1} +\dots + l_s d_{2r_s -1}$, which is contained in $\mathbb{Z} (p^\infty)$, can be represented in the form
\[
\phi_1 = \frac{z'}{p^{r_s^3 - r_s^2 -1}} + l_{s} \widetilde{f}_{r_s}, \mbox{ where } z' \in \mathbb{Z}.
\]
Since the number of summands in $\widetilde{f}_{r_s}$ is  $r_s +1 >10p(k+1)$ and $h-s <k+1$, there exists $r_s - 2> i_0 >2$ such that for every $1\leqslant w \leqslant h-s$ we have
\[
\mbox{either } r_{s+w} < r_s^3 - (i_0 +2)r_s  \mbox{ or }  r_{s+w} > r_s^3 - (i_0 -1)r_s .
\]
The set of all $w$ such that $r_{s+w} < r_s^3 - (i_0 +2)r_s$ we denote by $K$ (it may be empty or has the form $\{ 1,\dots ,a\}$ for some $1\leqslant a \leqslant h-s$). Set $L =\{ 1,\dots ,h-s \} \setminus K$. Thus
\[
\begin{split}
  \phi_0 = & \frac{z'}{p^{r_s^3 - r_s^2 -1}} + \sum_{w\in K} l_{s+w} d_{2r_{s+w}} + \frac{l_s}{p^{r_s^3 - r_s^2 }} +\dots +\frac{l_s}{p^{r_s^3 - (i_0 +2)r_s }} \\
   & +\left( \frac{l_s}{p^{r_s^3 - (i_0 +1)r_s}} +\frac{l_s}{p^{r_s^3 - i_0 r_s}}\right)+
    \left(\frac{l_s}{p^{r_s^3 - (i_0 -1)r_s}} + \dots +\frac{l_s}{p^{r_s^3 }} + \sum_{w\in L} l_{s+w} d_{2r_{s+w}} \right).
\end{split}
\]
So the element $\phi_2$  in row 1, can be represented in the form $\frac{c}{p^{r_s^3 - (i_0 +2)r_s}}$, for some $c \in \mathbb{Z}$.
Since $r_s > 10p(k+1)$, then $\frac{1}{1- 1/p^{r_s}} < \frac{32}{31}$ and $2k < p^{2k} <p^{r_s}$. Thus
we can estimate the second summand in row 2, which is denoted by $\phi_3$, as follows
\[
\begin{split}
|\phi_3 | = & \left| \left( \frac{l_s}{p^{r_s^3 - (i_0 -1)r_s}} + \dots +\frac{l_s}{p^{r_s^3 }} \right) + \sum_{w\in L} l_{s+w} \frac{1}{p^{r_{s+w}}} \right|  <
  \frac{|l_s |}{p^{r_s^3 - (i_0 -1)r_s}} \cdot \frac{1}{1-\frac{1}{p^{r_s}}} + \frac{k}{p^{r_s^3 - (i_0 -1)r_s +1} }\\
  < &  \frac{1}{p^{r_s^3 - (i_0 -1)r_s}} \left( k\frac{32}{31} + k\frac{1}{p} \right) <
\frac{2k}{p^{r_s^3 - (i_0 -1)r_s}} < \frac{1}{p^{r_s^3 - i_0 r_s}}.
\end{split}
\]
So
\[
\phi_0 = \frac{c}{p^{r_s^3 - (i_0 +2)r_s}}  +\frac{l_s}{p^{r_s^3 - (i_0 +1)r_s}} +\frac{l_s}{p^{r_s^3 - i_0 r_s}} + \phi_3 \not= 0
\]
and  if $\phi_0 =\frac{c''}{p^{\alpha}}, c'' \in \mathbb{Z}$, is an irreducible fraction, then $\alpha \geqslant r_s^3 - (i_0 +1)r_s > 5z$. Hence $\phi_0 \not= \frac{b}{p^z}$ and $\sigma \not= g$. Thus $\{ d_n\}$ is a $T$-sequence. $\Box$

Let us consider the group $\mathbb{Z} (p^{\infty})$ with discrete  topology. Then $\mathbb{Z} (p^{\infty})^\wedge =\Delta_p$ is the compact group of $p$-adic integers which elements are denoted by $x =(a_i) , 0\leqslant a_i <p$, and the identity is $\mathbf{1} =(1,0,0,\dots)$. By Remark 10.6 \cite{HR1}, $\langle \mathbf{1} \rangle$ is dense in $\Delta_p$ and, by 25.2 \cite{HR1}, $(\lambda , \mathbf{1})=\exp\{ 2\pi i \cdot \lambda \}$ for every $\lambda \in \mathbb{Z} (p^{\infty})$. We will use the following result, its proof see in Example 2.6.3 \cite{ZP2}. Let $d_{2n} =\frac{1}{p^n} \in \mathbb{Z} (p^{\infty})$ and $\widetilde{\mathbf{d}} =\{ d_{2n}\}$. Then $x \in s_{\widetilde{\mathbf{d}}} (\Delta_p)$ if and only if there exists $m=m(x) \in \mathbb{Z}$ such that
\begin{equation} \label{6}
(\lambda , x) = \exp ( 2\pi i m \lambda), \forall \lambda \in \mathbb{Z} (p^{\infty}).
\end{equation}
In other words, $x \in s_{\widetilde{\mathbf{d}}} (\Delta_p)$ iff $x = m\mathbf{1}$ for some $m\in\mathbb{Z}$. In particular, ${\rm Cl} \left(s_{\widetilde{\mathbf{d}}} (\mathbb{Z} (p^{\infty}))\right) = \Delta_p$.

{\bf Proof of Proposition \ref{p2}}.  The following two cases are fulfilled only.

 1) {\it $H_1= \mathbb{Z} (p^{\infty})\oplus H$}. Let $\mathbf{d} =\{ d_n\}$ be a $T$-sequence defined in Lemma \ref{l3}.
By Theorem B, it is enough to prove that $s_\mathbf{d} \left( (G_{\mathcal{D}})^\wedge \right) = \{ 0\}$.

Let
\[
\omega = x+y \in (G_{\mathcal{D}})^\wedge , \mbox{ where } x \in \mathbb{Z} (p^{\infty})^\wedge =\Delta_p ,  y\in (H_{\mathcal{D}})^\wedge  \mbox{ and } \; (d_n , \omega) \to 1.
\]
Then $(d_{2n} , \omega)= (\frac{1}{p^n},x)\to 1$. By (\ref{6}), $x = m \mathbf{1}$ for some $m\in \mathbb{Z}$  and $(\lambda, x) = \exp ( 2\pi i m \lambda), \forall \lambda \in \mathbb{Z} (p^{\infty})$. In particular, $(\widetilde{f}_n , x) = \exp ( 2\pi i m \widetilde{f}_n), \forall n\geqslant 3$. By (\ref{4}), we obtain that $(\widetilde{f}_n , x)\to 1$. Let $0\leqslant i< |H|-1$ be arbitrary. Then for any $n>i$ we have
\[
\begin{split}
\mbox{ if } |H|= \infty & :\quad
(d_{3(S_n+i)+1} , \omega)= (\widetilde{f}_{2j_{S_n +i}+1} , x)\cdot (e_i, y)\to (e_i, y) =1, \\
\mbox{ if } |H|< \infty & : \quad
(d_{3(n(|H|-1)+i)+1} , \omega)= (\widetilde{f}_{2j_{n(|H|-1) +i}+1} , x)\cdot (e_i, y)\to (e_i, y) =1.
\end{split}
\]
Since $e_i$ is arbitrary nonzero element of $H$, we have $y =0$.

Let us prove that also $x=0$. By hypothesis, $(d_{3n+2}, x) =(\beta_n , x)= \exp ( 2\pi i m \beta_n) \to \exp ( 2\pi i m \beta )=1$. Since $\beta$ is irrational, the equality $\exp ( 2\pi i m \beta )=1$ is possible only if $m=0$. Therefore $x=0$.

So $\omega =0$ and $s_\mathbf{d} \left( (G_{\mathcal{D}})^\wedge \right) = \{ 0\}$.

2) {\it $H_1= \langle e_0 \rangle \oplus H$, where $e_0$ is an (zero or nonzero) element of $\mathbb{Z} (p^{\infty})$}.

Let $\mathbf{d} =\{ d_n\}$ be a $T$-sequence defined in Lemma \ref{l4}.
By Theorem B, it is enough to prove that $s_\mathbf{d} \left( (G_{\mathcal{D}})^\wedge \right)^\perp = H_1 $.

Let
\[
\omega = x+y \in (G_{\mathcal{D}})^\wedge , \mbox{ where } x \in \mathbb{Z} (p^{\infty})^\wedge =\Delta_p ,  y\in (H_{\mathcal{D}})^\wedge  \mbox{ and } \; (d_n , \omega) \to 1.
\]
Then $(d_{2n} , \omega)= (\frac{1}{p^n},x)\to 1$. By (\ref{6}), $x = m \mathbf{1}$ for some $m\in \mathbb{Z}$  and $(\lambda, x) = \exp ( 2\pi i m \lambda), \forall \lambda \in \mathbb{Z} (p^{\infty})$. In particular, $(\widetilde{f}_n , x) = \exp ( 2\pi i m \widetilde{f}_n), \forall n\geqslant 1$. By (\ref{4}), we obtain that $(\widetilde{f}_n , x)\to 1$.

Let $0< i< |H|$ be arbitrary. Then for any $n>i$ we have
\[
\begin{split}
\mbox{ if } |H|= \infty & : \quad
(d_{2(S_n+i)-1} , \omega)= (\widetilde{f}_{S_n +i} , x)\cdot (e_i, y)\to (e_i, y) =1, \\
\mbox{ if } |H|< \infty & : \quad
(d_{2(n|H|+i)-1} , \omega)= (\widetilde{f}_{n|H| +i} , x)\cdot (e_i, y)\to (e_i, y) =1.
\end{split}
\]
Since $e_i$ is arbitrary nonzero element of $H$, $y =0$. Thus $\omega\in H^\perp$.

Let $i=0$. If $e_0 \not= 0$, then for any $n>i$ we have
\[
\begin{split}
\mbox{ if } |H|= \infty & : \quad
(d_{2S_n-1} , \omega)= (\widetilde{f}_{S_n} +e_0 , x)\to (e_0, x) =1, \\
\mbox{ if } |H|< \infty & : \quad
(d_{2n|H|-1} , \omega)= (\widetilde{f}_{n|H|} +e_0 , x) \to (e_0, x) =1.
\end{split}
\]
So $\omega\in \langle e_0\rangle^\perp$. Therefore $\omega \in H_1^\perp$.

If $e_0 =0$ and $H_1 =H$, we obtained that   $s_\mathbf{d} \left( (G_{\mathcal{D}})^\wedge \right) =\mathbb{Z} \mathbf{1}$ and it is dense in $\Delta_p =H^\perp$. Thus $s_\mathbf{d} \left( (G_{\mathcal{D}})^\wedge \right)^\perp = H_1 $.

Let $e_0 =\frac{c}{p^\delta} \not= 0$, where $c\in \mathbb{N}$. We need to prove only the converse inclusion, i.e. that $H_1^\perp \subset {\rm Cl} (s_{\mathbf{d}} ((G_{\mathcal{D}})^\wedge))$. Note that $H_1^\perp = \langle e_0 \rangle^{\perp} \cap (\mathbb{Z} (p^{\infty})^\wedge \oplus 0)$.
Set $x_0 =p^{\delta} \mathbf{1}$ and $y_0 =0$. Let $\omega = x_0 +y_0$. Then, by Remark 10.6 \cite{HR1}, $\langle\omega\rangle$ is dense in $H^\perp_1$. Hence it is enough to show that $\omega \in s_{\mathbf{d}} ((G_{\mathcal{D}})^\wedge)$. But this is indeed so since, by (\ref{6}), $(d_n, \omega) =(d_n, x_0)=\exp\{ 2\pi i \cdot p^{\delta} \lambda_n \}$, where $\lambda_{2n}= d_{2n}=\frac{1}{p^n}$ and
\[
\begin{split}
\mbox{ if } |H|= \infty & : \;
\lambda_{2S_n-1} = \widetilde{f}_{S_n} + e_0 = \widetilde{f}_{S_n} +\frac{c}{p^\delta} , \mbox{ and }  \lambda_{2(S_n+i)-1} =\widetilde{f}_{S_n +i} \mbox{ if } 0<i\leqslant n, \\
\mbox{ if } |H|< \infty & : \;
\lambda_{2n|H|-1} = \widetilde{f}_{n|H|} +e_0 =\widetilde{f}_{n|H|} + \frac{c}{p^\delta} , \mbox{ and }  \lambda_{2(n|H|+i)-1} =\widetilde{f}_{n|H| +i} \mbox{ if } 0<i<|H| -1.
\end{split}
\]
 Thus, by (\ref{4}), $(d_n, \omega ) \to 1$. $\Box$

The following lemma is a generalization of Lemma 2.4 \cite{Ga3}, in which the author modifies the construction of Example 5 \cite{ZP1} (or Example 2.6.2 \cite{ZP2}) and the proof of Proposition 3.3 \cite{Luk} (it is enough to put: $e'_i =e_i$ for $0\leqslant i<q$ and $e'_i =0$ for $i\geqslant q$).

\begin{lem} \label{l5}
Let $G= \bigoplus_{j=0}^\infty \langle e_j \rangle$, where $u_j := o(e_j) <\infty$ for every $j\geqslant 0$. Let $e'_i \in \langle e_i \rangle$ be arbitrary for each $0\leqslant i <\infty$.
Assume that one of the following conditions is fulfilled:
\begin{enumerate}
\item[a)] $u_{j_0} =u_{j_0 +1} =\dots$ for some $j_0 \geqslant 0$ and $u_{j_0}$ is divided on every $u_0,\dots, u_{j_0 -1}$. Set $b_k := e_k$.
\item[b)] $\sup\{ u_j, j\geqslant 0\} =\infty$. Let $e_{i_k}$ be such that $u_{i_k} \to \infty$. Put $b_k := e_{i_k}$.
\end{enumerate}
Set $\mu_n = S_{t(n)}$. Then the following sequence $\mathbf{d}=\{ d_n \} (n\geqslant 0)$
\[
\begin{split}
 d_{2n} : & \, e_0, 2 e_0, \dots, (u_0 -1) e_0, \; e_{1}, 2 e_{1}, \dots, (u_{1} -1) e_{1}, \dots \\ d_{1}=e'_0,  & \; d_3= e'_0 + b_1, \; d_5 = e'_1 +(b_2 + b_3), \\
 d_{2n+1} & = e'_{n\left({\rm mod}\, \mu_n \right) } + b_{S_{n-1} +1} +b_{S_{n-1} +2} + \dots + b_{S_{n}},
\end{split}
\]
is a $T$-sequence.
\end{lem}

\pr Let $k\geqslant 0$ and $g=\lambda_1 e_{v_1} + \lambda_2 e_{v_2} +\dots +\lambda_q e_{v_q} \not=0, v_1 <\dots <v_q$. We need to show that the condition of the Protasov-Zelenyuk criterion (Theorem 2.1.4 \cite{ZP2}) is fulfilled, i.e. there exists a natural $m$ such that $g\not\in  A(k,m)$. By the construction of $\mathbf{d}$, there is $m' >3$ such that $d_{2n} = \lambda(n) e_{r(n)}$, where $r(n) > \max \{ v_q , 3\}$ for every $n>m'$. Assume that $g\in A(k,2m_0)$ for some $m_0 >m'$. Then
\[
g=  l_1 d_{2r_1 +1} + l_2 d_{2r_2 +1} +\dots + l_s d_{2r_s +1} +
 l_{s+1} d_{2r_{s+1}} +l_{s+2} d_{2r_{s+2}} +\dots + l_{h} d_{2r_{h}},
\]
where all summands are nonzero, $\sum_{i=1}^h |l_i | \leqslant k+1$, $0<s\leqslant h$ (by our choosing of $m_0$) and
\[
\begin{split}
{} &  2m_0 < 2r_1 +1 < 2r_2 +1 < \dots < 2r_s +1 ,    \\
& 2m_0 \leqslant 2r_{s+1} < 2r_{s+2} < \dots < 2r_{h}.
\end{split}
\]
Moreover, by construction, all the elements $d_{2n+1} - e'_{n\left({\rm mod}\, \mu_n \right) }$ are independent. So, by the construction of the elements $d_{2n}$ and since $n\left({\rm mod}\, \mu_n \right) < S_{n-1}$ for every $n>3$, there is a subset $\Omega$ of the set $\{ s+1,\dots, h \}$ such that
\begin{equation} \label{7}
l_s (d_{2r_s +1} - e'_{r_s\left({\rm mod}\, \mu_n \right) }) + \sum_{w \in \Omega} l_w d_{2r_w} =l_s (b_{S_{r_s -1} +1} +b_{S_{r_s -1} +2} + \dots + b_{S_{r_s}}) +\sum_{w \in \Omega} l_w d_{2r_w}=0.
\end{equation}

a) {\it Assume that $u_{j_0} =u_{j_0 +1} =\dots$ and $u_{j_0}$ is divided on every $u_0,\dots, u_{j_0 -1}$}. Set $m_0 = 4m' (j_0 +2)(k+1)$. Then $d_{2r_s +1} - e'_{r_s\left({\rm mod}\, \mu_{r_s} \right) }$ contains exactly $r_s > m_0 -1 \geqslant 4k+4$
independent summands of the form $e_j$ with $j\geqslant S_{r_s -1} +1 > \frac{r_s(r_s-1)}{2} >m_0 > \max\{ m',j_0\}$.
Since $l_s d_{2r_s -1}\not= 0$ and $u_{j_0}$ is divided on every $u_0, \dots, u_{j_0 -1}$, we can assume that $l_s$ is not divided on $u_{j_0}$. So $l_s d_{2r_s -1}$ contains at least $4k+4$ non-zero independent  summands of the form $l_s b_j =l_s e_j, j> j_0$. Since $h-s \leqslant k$  and $l_w d_{2r_w}$ has the form $\lambda_w e_w$, the equality (\ref{7}) is impossible. Thus $g\not\in A(k,2m_0)$. So $\mathbf{d}$ is a $T$-sequence.

b) {\it Assume that $\sup\{ u_j, \; j\geqslant 0\} =\infty$}.  Choose $j' >m'$ such that $u_{i_j} > 2(k+1)$ for every $j> j'$. Set $m_0 = 4j'(k+1)$.
Then $d_{2r_s +1} - e'_{r_s\left({\rm mod}\, \mu_{r_s} \right) }$ contains exactly $r_s > m_0 -1 \geqslant 4k+4$ summands of the form $e_{i_j}$ with $j\geqslant S_{r_s -1} +1 > m_0 > j'$. So, since $|l_s| \leqslant k+1$, $l_s d_{2r_s -1}$ contains at least $4(k+1)$ non-zero independent  summands of the form $l_s e_{i_j}, j>j'$. Since $h-s \leqslant k$ and $l_w d_{2r_w}$ has the form $a_w e_w$, the equality (\ref{7}) is impossible. Thus $g\not\in A(k,2m_0)$. So $\mathbf{d}$ is a $T$-sequence. $\Box$

\begin{lem} \label{l6}
Let $G= \bigoplus_{j=0}^\infty \langle e_j \rangle$ endow with discrete topology and $o(e_j) <\infty$ for every $0\leqslant j <\infty$. Let $e'_j \in \langle e_j \rangle$ be arbitrary for every $0\leqslant j <\infty$. Set $H= \bigoplus_{j=0}^\infty \langle e'_j \rangle \subset G$ and $Y= \bigoplus_{j=0}^\infty \langle e'_j \rangle^\perp$, where $\langle e'_j \rangle^\perp$ is the annihilator of $\langle e'_j \rangle$ in $\langle e_j \rangle^\wedge$. Then $Y$ is dense in $H^\perp$.
\end{lem}

\pr It is clear that $Y \subset H^\perp$. Let $\omega = (a_0, a_1,\dots) \in H^\perp$. Then for any $j\geqslant 0$ and every $0\leqslant k < o(e'_j)$ we have $1=(\omega, k e'_j) = (a_j, k e'_j)$. Hence $a_j \in \langle e'_j \rangle^\perp$ for every $j\geqslant 0$. So $H^\perp \subset \prod_{j=0}^\infty \langle e'_j \rangle^\perp \subset \prod_{j=0}^\infty \langle e_j \rangle^\wedge =G^\wedge$. Since $Y$ is dense in $\prod_{j=0}^\infty \langle e'_j \rangle^\perp$, it is dense in $H^\perp$. $\Box$

{\bf Proof of Proposition \ref{p3}}. Let $\mathbf{d}=\{ d_n \}$ be the $T$-sequence which is defined in Lemma \ref{l5}. By Theorem B and Lemma \ref{6}, it is enough to show that
$s_{\mathbf{d}} ((G_{\mathcal{D}})^\wedge)) = \bigoplus_{j=0}^\infty \langle e'_j \rangle^\perp$.

We modify the proof of Proposition 3.3 \cite{Luk}.
Let $\omega =(a_0, a_1, \dots)\in s_{\mathbf{d}} ((G_{\mathcal{D}})^\wedge)$. By definition, there exists $N\in \mathbb{N}$ such that $|1- (d_{2n}, \omega)|<0.1, \forall n>N$. Thus, there is $N_0 >N$ such that $|1- ( l e_{j}, \omega)|= |1- (l e_{j}, a_j)| <0.1, \forall l=1,\dots, u_j -1,$ for every $j>N_0$. This means that $a_j =0$ for every $j>N_0$. So $\omega \in \bigoplus_{j=0}^\infty \langle e_j \rangle \subset (G_{\mathcal{D}})^\wedge$. Set $\mu_n = S_{t(n)}$. Since  $(d_{2(\mu_n +j)-1}, \omega) \to 1$ at $n\to\infty$ too and $(d_{2(\mu_n +j)-1}, \omega) = (e'_j, a_j)$ for all sufficiently large $n$, we obtain that $a_j \in \langle e'_j \rangle^\perp$ for any $j\geqslant 0$. Thus $ s_{\mathbf{d}} ((G_{\mathcal{D}})^\wedge) \subset \bigoplus_{j=0}^\infty \langle e'_j \rangle^\perp$. The converse inclusion is trivial. Hence $ s_{\mathbf{d}} ((G_{\mathcal{D}})^\wedge) = \bigoplus_{j=0}^\infty \langle e'_j \rangle^\perp$. $\Box$

{\bf Proof of Theorem \ref{t2}}. {\it Let us prove the implication $(1)\Rightarrow (2)$.}

Since $G$ is bounded, it is a finite sum of its $p$-subgroups $G_p$. So we may assume that $G=\bigoplus_{p\in M} G_p$ and $H=\bigoplus_{p\in M} H_p$, where $M$ is a finite set of prime numbers.
Let $N=\exp (H) =p_1^{b_1}  \dots  p_l^{b_l}$, where $p_1, \dots, p_l$ be distinct prime integers.
Then $G$ contains a subgroup of the form $\mathbb{Z} (N)^{(\omega)}$ if and only if $G_{p_i}$ contains
$\mathbb{Z} (p_i^{b_i})^{(\omega)}$ for every $1\leqslant i\leqslant l$.
Since the von Neumann radical of a finite product of topological groups is the product of their von Neumann radicals, we may assume that $G$ is a $p$-group. Let us note that $G$ contains a subgroup of the form $\mathbb{Z} (p^n)^{(\omega)}$ iff $G$ contains infinitely many independent elements of order greater or equal than $p^n$ (Lemma 8.1 \cite{Fuc}).

By the  Pr\"{u}fer Theorem (Theorem 11.2 \cite{Fuc}), we may assume that $H$ has the form $H =\bigoplus_{i=1}^{n} \mathbb{Z}(p^{a_i})^{(k_{i})}$, where $k_i >0$ and $1\leqslant a_1 <\dots < a_n$.

Let $e_{j,i}, 0\leqslant j < k_i, 1\leqslant i\leqslant n,$ be such that
\[
H =\bigoplus_{i=1}^{n} \bigoplus_{0\leqslant j< k_i} \langle e_{j,i} \rangle, \mbox{ and } o(e_{j,i}) = p^{a_i}, 1\leqslant i\leqslant n.
\]
Let  $B_1$ be the set of all indexes $i$ such that $k_{i}<\infty$ and $B_2$ be the set of all  $i$ for which $k_{i}=\infty$. Set $q=\sum_{i\in B_1} k_i$.

1) Let $k_n =\infty$. Then
\[
H =\left( \bigoplus_{i\in B_1} \mathbb{Z}(p^{a_i})^{(k_{i})} \oplus \mathbb{Z}(p^{a_n})^{(\omega)} \right) \oplus \bigoplus_{i\in B_2, i\not= n} \mathbb{Z}(p^{a_i})^{(\omega)}
\]
and it is a finite sum of  groups of the form a) in Proposition \ref{p3}. So $H$ admits a complete MinAP  group topology.

2) Let $k_n <\infty$. Denote by $b_{j,i}$ the height of $e_{j,i}$ in $G$. Set $C_i$ is the set of all indexes $j$ such that $b_{j,i}\geqslant p^{a_n}$ (it may be empty) and put $C=\cup_{i=1}^n C_i$. Then the following cases can be fulfilled only.

a) {\it The cardinality $|C|$ of $C$ is finite}. By hypothesis, there is a sequence $\{ e_m \}_{m=q}^\infty$ such that $o(e_m)= k_n$ and all elements $e_m$ and $e_{j,i}$ are independent. Set
\[
 G_1=\left( \bigoplus_{i\in B_1} \bigoplus_{0\leqslant j< k_i} \langle e_{j,i} \rangle \oplus \bigoplus_{m=q}^\infty \langle e_m\rangle \right) \oplus \bigoplus_{i\in B_2} \mathbb{Z}(p^{a_i})^{(\omega)} .
\]
By Reduction Principle and Proposition \ref{p3}(a), there is  a complete Hausdorff  group topology $\tau$ on $G$ such that $\mathbf{n}(G,\tau) =H$.

b) $|C|=\infty$. Then there is an $i_0\in B_2$ such that $|C_{i_0} |=\infty$. Let us denote by $e_j$ an element of order $p^{a_n}$ such that $p^{a_n - a_{i_0}} e_j = e_{j,i_0}, j\in C_{i_0}$. Set $D_{i_0} = \mathbb{N}\setminus C_{i_0}$. Then all elements $e_m$ and $e_{j,i},$ where either $i\not= i_0$ or $i=i_0$ and $j \in D_{i_0},$ are independent.

If $|D_{i_0} |=\infty$, we set
\[
G_1 = \left( \bigoplus_{i\in B_1} \bigoplus_{0\leqslant j< k_i} \langle e_{j,i} \rangle \oplus \bigoplus_{j\in C_{i_0}} \langle e_j \rangle \right) \oplus \bigoplus_{i\in B_2, i\not= i_0} \mathbb{Z}(p^{a_i})^{(\omega)} \oplus \bigoplus_{j\in D_{i_0}} \langle e_{j,i_0} \rangle .
\]
For the first group we can put $e'_l = e_{j,i}$ if $i\in B_1, 0\leqslant j <k_i, 0\leqslant l < q$, and $e'_l = e_{j,i_0}$ if $j\in C_{i_0}$ and $l\geqslant q$. By Reduction Principle and Proposition \ref{p3}(a), we can find  a complete Hausdorff  group topology $\tau$ on $G$ such that $\mathbf{n}(G,\tau) =H$.

If $|D_{i_0} |<\infty$, we put
\[
G_1 = \left( \bigoplus_{i\in B_1} \bigoplus_{0\leqslant j< k_i} \langle e_{j,i} \rangle \oplus \bigoplus_{j\in D_{i_0}} \langle e_{j,i_0} \rangle \oplus \bigoplus_{j\in C_{i_0}} \langle e_j \rangle \right) \oplus \bigoplus_{i\in B_2, i\not= i_0} \mathbb{Z}(p^{a_i})^{(\omega)}  .
\]
Analogously, by Reduction Principle and Proposition \ref{p3}(a), there is  a complete Hausdorff  group topology $\tau$ on $G$ such that $\mathbf{n}(G,\tau) =H$.

{\it The implication $(2)\Rightarrow (3)$ is trivial.}

{\it Let us prove the implication $(3)\Rightarrow (1)$.} Let
\[
\exp H =p_1^{b_1}  \dots  p_l^{b_l} \mbox{ and } \exp G =p_1^{n_1}  \dots  p_l^{n_l} \cdot p_{l+1}^{n_{l+1}}  \dots  p_t^{n_t},
\]
where $p_1, \dots, p_t$ be distinct prime integers.

Assuming the converse, we obtain that there is a $1\leqslant j \leqslant l$ such that $G$ contains only a finite subset of independent elements $g$ for which $o(g) = p_j^a$ with $a\geqslant b_j$. Set $m := \exp G /p_{j}^{n_j - b_j+1}$ and $\pi : G\to G, \, \pi (g) = m g$. Then $\pi (H) \not= 0$ and, by our hypotheses, $\pi (G)$ is finite.

Let us prove that there is no  Hausdorff group topology $\tau$ such that $\mathbf{n}(G,\tau)=H$. (We repeat the arguments of D.~Remus (see \cite{Com}).

Let $\tau$ be any Hausdorff group topology on $G$. Then ${\rm Ker} (\pi)$ is open and closed and, hence, dually closed and dually embedded \cite{Nob}. So, by Lemma \ref{l1}, $\mathbf{n} (G,\tau) \subset {\rm Ker} (\pi)$. Hence $H\not= \mathbf{n} (G,\tau)$.
This completes the proof of the theorem. $\Box$

{\bf Proof of Corollary \ref{c4}}. By Reduction Principle we may assume that $G= \langle e_0 \rangle + H$. Since $\langle e_0 \rangle \cap H=\{ 0\}$ and $e_0$ of infinite order, $G=\mathbb{Z} \oplus H$. So the assertion follows from Proposition \ref{p1}. $\Box$

{\bf Proof of Corollary \ref{c5}}. By Baer's Theorem and Reduction Principle we may assume that $G= \mathbb{Z} (p^\infty) + H= \mathbb{Z} (p^\infty) \oplus H_1$, where $H\supset H_1$. So the assertion follows from Proposition \ref{p2}. $\Box$

{\bf Proof of Proposition \ref{p4}}. Let us endow the group $G$ under the asterisk group topology \cite{Kap}. By the Kaplan Theorem \cite{Kap}, $G^\wedge = \prod_{\alpha \in I} G_\alpha^\wedge$. Since $G_\alpha^\wedge$ is trivial, then also  $G^\wedge$ is trivial. $\Box$

{\bf Proof of Theorem \ref{t3}}.
By the  Pr\"{u}fer Theorem (Theorem 11.2 \cite{Fuc}), we may assume that $G$ has the form $G =\bigoplus_{p\in M} \bigoplus_{i=1}^{n_p} \mathbb{Z}(p^{a_{i,p}})^{(k_{i,p})}$, where $k_{i,p} >0$ and $1\leqslant a_{1,p} <\dots < a_{n_p,p}$ and $M$ is a finite set of prime numbers.

{\it Let us prove that $G$ admits a MinAP group topology}. Since the von Neumann radical of a finite product of topological groups is the product of their von Neumann radicals, we may assume that $G$ is a $p$-group. So we may assume that $G$ has the form
\[
G =\bigoplus_{i=1}^{n} \mathbb{Z}(p^{a_i})^{(k_{i})}, \mbox{ where } k_i >0 \mbox{ and } 1\leqslant a_1 <\dots < a_n .
 \]
Let  $B_1$ be the set of all indexes $i$ such that $k_{i}<\infty$ and $B_2$ be the set of all  $i$ for which $k_{i}=\infty$. Since, by hypothesis, $k_n =\infty$, we can represent $G$ in the form
\[
 G=\left( \bigoplus_{i\in B_1} \mathbb{Z}(p^{a_i})^{(k_{i})} \oplus \mathbb{Z}(p^{a_n})^{(\omega)} \right) \oplus \bigoplus_{i\in B_2}  \bigoplus_{\widetilde{k}_{i}} \mathbb{Z}(p^{a_i})^{(\omega)},
\]
where the cardinal number $\widetilde{k}_{i}$ is infinite and such that $k_i = \omega \cdot \widetilde{k}_{i}$ for every $i\in B_2$. By Corollary \ref{c2}, each group in this representation of $G$ admits a complete MinAP group topology. Thus, by Proposition \ref{p4}, $G$ admits a MinAP group topology.

{\it Let us prove the converse assertion.} (We repeat the arguments of D.~Remus (see \cite{Com}). Assuming the converse, we obtain that there is  $p_0 \in M$ such that $k_{n_{p_0}, p_0}< \infty$. Set $m := \exp G /p_{0}^{a_{n_{p_0}, p_0} -1}$ and $\pi : G\to G, \, \pi (g) = m g$. Then $\pi (G) \not= 0$ and, by our hypotheses, $\pi (G)$ is finite.

Now let $\tau$ be any Hausdorff group topology on $G$. Then ${\rm Ker} (\pi)$ is open and closed and, hence, dually closed and dually embedded \cite{Nob}. So, by Lemma \ref{l1}, $\mathbf{n} (G,\tau) \subset {\rm Ker} (\pi)$. Hence $G\not= \mathbf{n} (G,\tau)$. It is a contradiction.
This completes the proof of the theorem. $\Box$

\bibliographystyle{plain}

\end{document}